\documentclass[12pt,reqno]{amsart}
\usepackage{amssymb,delarray}
\usepackage{amsfonts}
\usepackage{epsfig}
\usepackage[all]{xy}

\def\leq{\leqslant}
\def\geq{\geqslant}

\def\l{\operatorname{l}}

\newtheorem{thm}{Theorem}
\newtheorem{lem}
{Lemma}
\newtheorem{prop}
{Proposition}
{Claim}

\newtheorem{cor}
{Corollary}
\newtheorem{rem}
{Remark}
{Question}

{\catcode`\@=11
\gdef\n@te#1#2{\leavevmode\vadjust{%
 {\setbox\z@\hbox to\z@{\strut#1}%
  \setbox\z@\hbox{\raise\dp\strutbox\box\z@}\ht\z@=\z@\dp\z@=\z@%
  #2\box\z@}}}
\gdef\leftnote#1{\n@te{\hss#1\quad}{}}
\gdef\rightnote#1{\n@te{\quad\kern-\leftskip#1\hss}{\moveright\hsize}}
\gdef\?{\FN@\qumark}
\gdef\qumark{\ifx\next"\DN@"##1"{\leftnote{\rm##1}}\else
 \DN@{\leftnote{\rm??}}\fi{\rm??}\next@}}

\begin{document}
\baselineskip=13.7pt plus 2pt 

\title[Numerically
pluricanonical coverings] {On numerically pluricanonical cyclic
coverings}
\author[V. Kharlamov and Vik.S. Kulikov ]{V. Kharlamov and Vik.S. Kulikov }

 \address{Strasbourg University}
\email{ kharlam@math.unistra.fr}

\address{Steklov Mathematical Institute}
 \email{kulikov@mi.ras.ru}

\dedicatory{} \subjclass{}
\thanks{ The first author was partially funded by the
ANR-09-BLAN-0039-01 grant of Agence Nationale de la Recherche. The
second author was partially supported by grants of NSh-4713.2010.1,
RFBR 11-01-00185, and by AG Laboratory HSE, RF government grant, ag.
11.G34.31.0023.}


\begin{abstract}
In this article, we investigate some properties of cyclic coverings
$f: Y\to X$ of complex surfaces of general type $X$ branched along
smooth curves $B\subset X$ that are numerically equivalent to a
multiple of the canonical class of $X$. The main results concern
coverings of surfaces of general type with $p_g=0$ and Miyaoka--Yau
surfaces; in particular, they provide new examples of
multicomponent  moduli spaces of  surfaces with given Chern numbers
as well as new examples of surfaces that are not deformation
equivalent to their complex conjugates.

\end{abstract}

\maketitle
\setcounter{tocdepth}{1}


\def\st{{\sf st}}



\section*{Introduction} \label{introduc}

In this article, we investigate some properties of cyclic coverings
of algebraic surfaces. Let us recall that there exist two equivalent
approaches to define and to treat branched coverings. One is due to
Grauert-Remmert-Stein theorem, which insures that the covering is
uniquely defined by its unramified part; in particular, if the base
of a branched covering is nonsingular, then the branch locus is
either empty or a divisor, and each unbranched finite covering of
the complement of a divisor extends uniquely to a branched covering
with a normal covering variety. The other, also traditional, one, to
which we give preference in this paper, is due to the canonical
equivalence between finite branched coverings of a given nonsingular
(or normal) variety $X$ and the finite field extensions of its
rational function field $\mathbb C(X)$. A branched covering is
called {\it Galois}, if the field extension is Galois, and it is
called {\it cyclic} if the Galois group is cyclic.

Thus, given a finite cyclic cover $Y\to X$, one can speak on the
action of a finite cyclic group $G$ on $Y$ and identify $X$ with the
quotient variety $X=Y/G$. It is worth, however, to underline that
when we speak on a cyclic Galois group and a cyclic Galois covering,
we are not fixing an isomorphism between the Galois group and a
group of permutations of a given finite set; in particular, we call
two Galois coverings, $f_1: Y_1\to X$ and $f_2:Y_2\to X$, isomorphic
if, and only if, there exists an isomorphism $\phi: Y_1\to Y_2$ such
that $f_2\circ \phi=f_1$ and $\phi$ transforms Galois automorphisms
in Galois automorphisms.

The cyclic coverings that we consider in this paper are somehow
special. Namely, they have a non-empty branch locus $B\subset X$ and
we forbid the Galois group action to have points in $Y$ whose
stabilizer is a non-trivial proper subgroup of the Galois group. We
call such coverings {\it totally ramified}. (Note, that since we
consider only the coverings over a smooth $X$, the above restriction
forbids also the Galois group to have isolated fixed points in $Y$.)

Note that the assumption to be totally ramified being applied to the
points of the ramification divisor $R\subset Y$ implies that
$f^*(B)=dR$, $d=\deg f$. Furthermore, a degree $d$ cyclic covering
$f:Y\to X$ is totally ramified (over $X$) if, and only if,
$f^*(B)=dR$ and $f$ is unramified over $X\setminus B$. Note that if
$X_1$ is birationally equivalent to $X$ then the cyclic covering
$f_1:Y_1\to X_1$,
induced by a totally ramified 
cyclic covering $f:Y\to X$, is not necessary totally ramified over
$X_1$ if the branch curve $B\subset X$ of $f$ has singular points.

By a {\it numerically multi-canonical} cyclic covering we understand
a totally ramified cyclic covering whose branch curve is
non-singular and numerically equivalent to a multiple of the
canonical class.

The main results of the article can be divided into three groups.

First, we consider the moduli space of surfaces with given square of
the canonical class, $K^2_Y$, and given arithmetic genus, $p_a(Y)$,
and, under assumption that it contains surfaces given by $d$-sheeted
totaly ramified numerically multi-canonical cyclic coverings $f:Y\to
X$ of a surface of general type $X$, provide a lower bound for the
number of connected components of such a moduli space depending on
the number of elements of the torsion group $\text{Tor}H^2(X,\mathbb
Z)$ (see Theorem \ref{mod}, 
Proposition \ref{Camp1}, Remark \ref{rem+}, and Corollaries
\ref{Cnew}, \ref{Cplus}).

Second,  we investigate  the degree of the canonical map
$\varphi_{K_Y}:Y\to \mathbb P^{p_g(Y)-1}$ for surfaces $Y$ given as
two-sheeted totaly ramified numerically multi-canonical cyclic
coverings $f:Y\to X$, where $X$ is a surface of general type with
$p_g(X)=0$ (see Theorem \ref{pg0}, Corollary \ref{pg01}, and
Propositions \ref{Camp} -- \ref{Me-Pa}).

Third, we show that if a Miyaoka--Yau surface $X$ (that is, a
surface of general type with $c_1^2=3c_2$) has no anti-holomorphic
automorphisms (such a property is shared, for example, by all  fake
projective planes, see \cite{KK1}) and if a surface $Z$ is
deformation equivalent to $Y$ given as a totaly ramified numerically
multi-canonical cyclic covering $f:Y\to X$, then $Z$ also has no
anti-holomorphic automorphisms (see Theorem \ref{main1}). Note that,
together with the known examples
(see, \cite{Ma}, \cite{KK1}, \cite{KK2}, \cite{Ca}) of pairs of
complex surfaces $(Z,Z')$ that are not deformation equivalent
despite being diffeomoprhic (with preserving the orientation), the
surfaces $Y$ involved in Theorem \ref{main1} and their complex
conjugates, $\bar Y$, give infinite series of new examples of
diffeomoprhic, but not deformation equivalent, surfaces. Besides, we
prove that those connected components $M$ of the moduli space of
surfaces that contain surfaces $Y$, obtained as two-sheeted totaly
ramified numerically multi-canonical cyclic coverings $f:Y\to X$ of
a Miyaoka-Yau surface $X$ branched along a curve $B$ numerically
equivalent to $2mK_X$, are irreducible varieties of complex
dimension $\dim M= m(2m-1)K_X^2+p_g(X)$ and Kodaira dimension
$\kappa (M)=-\infty$, moreover, if the irregularity $q(X)=0$ then
$M$ is an unirational variety (see Theorem \ref{main3} and also
Remark \ref{last-rem}).

Up to ou knowledge, there are a very few, referred to in the
literature, examples of complex surfaces of general type with such a
nontrivial action of a group that is deforming simultaneously with
any deformation of complex structure (cf., \cite{KK2}). As a
consequence of the proof of Theorem \ref{main3} we provide infinite
series of such surfaces with action of the group $\mathbb Z/2\mathbb
Z$ (see Remark \ref{last-rem}).

\section{Totally ramified cyclic coverings: fundamental groups and classification}

Recall that each continuous map $f:V\to W$ of path connected
topological spaces defines the homomorphism of fundamental groups,
$f_* : \pi_1(V,q)\to \pi_1(W,p)$, for any pair of points, $q\in V$
and $p\in W$ with $p=f(q)$, and that for any pair of basic points,
$(q',p')$ and $(q'',p'')$ connected by the pathes, $h$ in $V$ and
$f(h)$ in $W$, these homomorphisms are conjugated by
change-of-basepoint homomoprhisms $\beta_h$ and $\beta_{f(h)}$ in
the sense that $f^{\prime\prime}_{*}\circ \beta_h=\beta_{f(h)}\circ
f'_{*}$. Note that all these homomorphisms $f_*$ induce one and the
same homomorphism
$$H_1(V,\mathbb Z)\simeq \pi_1(V,q)/[\pi_1(V,q),\pi_1(V,q)]\to
H_1(W,\mathbb Z)\simeq \pi_1(W,p)/[\pi_1(W,p),\pi_1(W,p)]$$ at the
level of the first homology groups. Below, when the choice of the
base points $p\in W$ and $q\in f^{-1}(p)$ in the fundamental groups
$\pi_1(V,q)$ and $\pi_1(W,p)$ is not essential for the proof, we
omit mentioning the base points and denote the fundamental groups by
$\pi_1(V)$ and $\pi_1(W)$, respectively.

\begin{prop} \label{pi}  Let $X$ be an irreducible smooth projective surface, let
$B\subset X$ be an irreducible reduced smooth curve divisible by $d$
as an element of $\text{Pic}(X)$,  and let $f:Y\to X$ be a
$d$-sheeted totally ramified cyclic covering branched along $B$. If
$(B^2)_X>0$ then the homomorphism $f_*: \pi_1(Y)\to \pi_1(X)$
induced by  $f$ is an isomorphism.
\end{prop}

\begin{lem}\label{NoriLemma}
Let $B$ be an irreducible reduced smooth curve on an irreducible
smooth projective surface $X$. If $(B^2)_X>0$, then the kernel $K$
of the epimorphism $\pi_1(X\setminus B)\to \pi_1(X)$ is a cyclic
group and it is a subgroup of the center of $\pi_1(X\setminus B)$.
\end{lem}
\proof  Let $N\subset X$ be the closure of a tubular neighbourhood
of $B$. Then $N$ has a structure of a locally trivial
$C^{\infty}$-fibration over $B$ with fiber $D=\{ z\in \mathbb C\mid
|z|\leq 1\}$. If we delete one fiber of this fibration (say, over a
point $b\in B$), then we get a trivial fibration $N_0\simeq D\times
B_0$ over $B_0=B\setminus\{b\}$. Consider the section $B_1=B_0\times
\{ z=1\}$ of the latter fibration and choose a point $p\in B_1$.

Let $\gamma\in G=\pi_1(X\setminus B,p)$ be an element represented by
the loop $\partial D$, where $\partial D$ is the boundary of the
fibre of $N\to B$ containing $p$. We have the exact sequence of
groups
\begin{equation}\label{ex1}
1\to K\to G\to \pi_1(X,p)\to 1,
\end{equation}
where $K$ is the normal closure of $\gamma$ in $G$.

By Nori's Weak Lefschetz Theorem (see \cite{No}, Proposition 3.27),
$K$ is a finitely generated abelian group and its centralizer $C(K)$
in $G$ is a subgroup of finite index.

On the other hand, we have $\pi_1(N_0\setminus
B_0,p)=\pi_1(B_1,p)\times \{\gamma^n \mid n\in \mathbb Z\}$. Again
by Nori's Weak Lefschetz Theorem (see \cite{No}, Proposition 3.26),
the embedding $i:B\hookrightarrow X$ induces an epimorphism
$i_*:\pi_1(B)\to\pi_1(X)$. Therefore, the elements of $\pi_1(B_1,p)$
generate the group $\pi_1(X,p)$ and the embedding $N_0\setminus
B_0\hookrightarrow X\setminus B$ induces an epimorphism
$$\pi_1(N_0\setminus B_0,p)=\pi_1(B_1,p)\times \{\gamma^n \mid n\in \mathbb Z\}\to
\pi_1(X\setminus B,p).$$ Hence, the normal closure $K$ of $\gamma$
in $\pi_1(X\setminus B,p)$ is nothing but the cyclic group generated
by
$\gamma$ and, furthermore, $K$ is contained in the center of $\pi_1(X\setminus B,p)$. \qed \\

\proof[Proof of Proposition \ref{pi}] The covering $f:Y\to X$
induces an embedding of the fundamental groups $f_*:\pi_1(Y\setminus
f^{-1}(B))\hookrightarrow \pi_1(X\setminus B)$ such that $\overline
G=f_*(\pi_1(Y\setminus f^{-1}(B)))$ is a subgroup of
$G=\pi_1(X\setminus B)$ of index $d=\deg f$.

Let $\gamma$ denote, as in the proof of Lemma \ref{NoriLemma}, a
generator of $K$. Since $f^*(B)=dR$, we have  $\gamma^i\not \in
\overline G$ for $1\leq i\leq d-1$ and $\gamma^d\in \overline G$.
Therefore, $G=\overline G\cup \gamma\overline G\cup \dots \cup
\gamma^{d-1}\overline G$. On the other hand, $\pi_1(Y)\simeq
\overline G/\langle \gamma^d\rangle$. Exact sequence (\ref{ex1})
gives rise to the exact sequence
\begin{equation}\label{ex2}
1\to K/\langle \gamma^d\rangle\to G/\langle \gamma^d\rangle\to
\pi_1(X)\to 1.
\end{equation}
The group $\overline G/\langle \gamma^d\rangle \simeq \pi_1(Y)$
naturally imbeds into $G/\langle \gamma^d\rangle$, so that
Proposition \ref{pi} follows now from exact sequence (\ref{ex2}) and
the equality
$$G/\langle
\gamma ^d\rangle=\pi_1(Y)\cup \gamma\pi_1(Y)\cup \dots \cup
\gamma^{d-1}\pi_1(Y). \qquad \qquad \qed $$

Let $C_1$ and $C_2$ be two divisors on a surface $X$. We will use
the notation $C_1\sim C_2$ if $C_1$ and $C_2$ are linear equivalent
and $C_1\equiv C_2$ if $C_1$ and $C_2$ are numerically equivalent.
Denote by $\text{Tor}_d\text{Pic}(X)$ the subgroup of the Picard
group $\text{Pic}(X)$ consisting of the elements whose order is a
divisor of $d$, and by $N_{X,d}=|\text{Tor}_d\text{Pic}(X)|$ the
order of $\text{Tor}_d\text{Pic}(X)$. Given a divisor $B$ whose
divisor class is divisible by $d\in\mathbb N$, we denote by
$(B)_d\subset \text{Pic}(X)$ the set of divisor classes $\beta$ such
that $B\sim d\beta$. Clearly, $(B)_d$ is a principal homogeneous
space over $\text{Tor}_d\text{Pic}(X)$.

The following lemma is well-known. Since we did not find an
appropriate reference, we supply this lemma with a proof.
\begin{lem}
\label{cov} Let $B\subset X$ be an irreducible  reduced curve. If
$B$, as an element of $\text{Pic}(X)$, is divisible by $d$ {\rm
(}that is $B\sim dC$ for some divisor $C${\rm )}, then there is a
natural bijection between $(B)_d$ and the set of isomorphism classes
of $d$-sheeted totally ramified cyclic coverings $f_i:Y_i\to X$
branched along $B$. If $[C_i]\in (B)_d$ is the divisor class
corresponding to $f_i$ under this bijection, then $[R_i]=f_i^*[C_i]$
and  $d[C_i]=[B]$.

If $B$ is a smooth curve, then each of these $Y_i$ is a smooth
surface.
\end{lem}
\proof By definition,  a finite morphism $f:Y\to X$ is a branched
$d$-sheeted cyclic covering of a smooth surface $X$ if $Y$ is a
normal surface and $f^*: \mathbb C(X)\hookrightarrow \mathbb C(Y)$
is a finite Galois extension of fields with Galois group
$\text{Gal}(Y/X)\simeq \mathbb Z/d\mathbb Z$. Let $h^*$ be a
generator of $\text{Gal}(Y/X)$. In accordance with so called
Hilbert's theorem 90, in the field $\mathbb C(Y)$,
considered as the vector space over $\mathbb C(X)$, one can choose a
basis $w_0=1,w_1,\dots,w_{d-1}$ over $\mathbb C(X)$ such that
$h^*(w_i)= \mu^i w_i$, where $\mu$ is a primitive $d$-root of unity.
Hence, we can put $w=w_1$ and conclude that each branched
$d$-sheeted cyclic covering can be seen as an extension $\mathbb
C(Y)=\mathbb C(X)(w)$ with $h^*(w)= \mu w$, and $w^d= g\in \mathbb
C(X)$ for some function $g$. The reverse statement is also
straightforward.

The automorphism $h^*$ defines an automorphism $h:Y\to Y$, while the
branch locus of $f$ consists of those irreducible curves $D_i$ that
participate in the principal divisor  $(g) =\sum_{i=0}^n a_iD_i \in
\text{Div}(X)$, $a_i\in \mathbb Z$, with coefficients $a_i\not
\equiv 0\,\text{mod}\, d$. Furthermore, $gcd(a_i,d)$ is equal to the
number of points in the inverse image of a generic point in $D_i$.
Therefore, in our case there is a unique curve $D_i$ (say, $D_0$;
and it must be the branch curve $B$) for which $a_0\not \equiv
0\,\text{mod}\, d$ and, moreover, $a_0$ and $d$ are coprime. Thus,
we can find an integer $b$ coprime with $d$ such that $(g^b)=B-dC$
for some divisor $C\in \text{Div}(X)$. Therefore, choosing a
function $\tilde g\in \mathbb C(X)$ and replacing $w$ by $w^b \tilde
g$, $g$ by $\tilde g^d g^b$, and $\mu$ by $\mu^b$, we get a
presentation $\mathbb C(Y)=\mathbb C(X)(w)$, $h^*(w)=\mu w$, $w^d=g$
with $(g)=B-dC'$, where $C'$ is any given divisor linear equivalent
to $C$. This provides a natural bijection between $(B)_d$ and the
set of isomorphism classes of $d$-sheeted cyclic coverings of $X$
branched along $B$.

If $R$ is the ramification locus of $f$, $R=f^{-1}(B)$ and
$f^*(B)=dR$ for an irreducible curve $B$, then
$(w^d)=(f^*(g))=dR-df^*(C)$. Therefore, $(w)=R- f^*(C)$. \qed

\begin{rem} \label{Remark} {\rm To enumerate the branched coverings as in
Lemma \ref{cov}, it is often convenient to choose (as in the proof
of Lemma \ref{cov}) as "a base point" some divisor $C$ such that
$B\sim dC$ and then make use of the bijection that associates with
$\alpha_i\in \text{Tor}_d\text{Pic}X$ the $d$-sheeted cyclic
coverings defined by $w^d=g_i$, where} $(g_i)=B-dC_i$ with $C_i\sim
C+\alpha_i$. \end{rem}

\begin{rem}\label{rem3} {\rm Lemma \ref{cov}
and Remark \ref{Remark} hold also if $B$ is not necessary
irreducible, but reduced, curve and $d=2$.}
\end{rem}

Let $f:Y\to X$ be a $d$-sheeted totally ramified cyclic covering
$f:Y\to X$ branched along a smooth irreducible curve $B\subset X$
with $(B^2)_X>0$. By Proposition \ref{pi}, the map $f$ induces an
isomorphism between the fundamental groups, $\pi_1(Y)$ and
$\pi_1(X)$. Therefore, it induces also an isomorphism $f_*$ between
$H_1(Y,\mathbb Z)$ and $H_1(X,\mathbb Z)$ and, in particular,
between their torsion subgroups, $\text{Tor} H_1(Y,\mathbb Z)$ and
$\text{Tor} H_1(X,\mathbb Z)$, which in its turn implies, by the
universal coefficient theorem (that gives $\text{Tor}H_1(V,\mathbb
Z)=\text{Ext}(H_1(V;\mathbb Z),\mathbb Z)=\text{Tor}H^2(V,\mathbb
Z)$ for every compact manifold $V$) and its functoriality, that $f^*
: \text{Tor}H^2(X,\mathbb Z)\to \text{Tor}H^2(Y,\mathbb Z)$ is an
isomorphism as well.

In fact, even a stronger statement holds under our assumptions.
Indeed, the exponential exact sequence of sheaves
$$0\to \mathbb Z\to \mathcal O_Y\to \mathcal O_Y^*\to 0$$
provides a short exact sequence
$$
0\to \text{Tor}\text{Pic}^0(Y)\to \text{Tor Pic}(Y)
\xrightarrow{\delta} \text{Tor}H^2(Y,\mathbb Z)\to 0,
$$
where $\text{Pic}^0(Y)$ is the connected component of $0$ in $\text{
Pic}(Y)$ and $\delta$ is the connecting (first Chern class)
homomorphism. Hence, using the functorial isomorphism
$\text{Pic}^0(V)=H^{0,1}(V)/(H^1(V;\mathbb Z)^{0,1}$ and the five
lemma applied to the diagram  formed by the above short exact
sequence and its copy written for $X$, one proves that $f^* :
\text{Tor Pic }(X)\to \text{Tor Pic}(Y)$ is also an isomorphism.

In what follows,  we use only $\text{Tor}H^2(V,\mathbb Z)$ and, for
shortness, abbreviate this  notation to $\text{Tor}(V)$.

\section{Numerically multi-canonical cyclic coverings}
In this section we fix an integer $d\geq 2$ and a smooth irreducible
curve $B\equiv dmK_X$ on a nonsingular surface $X$. A $d$-sheeted
totally ramified cyclic covering $f:Y\to X$ branched along $B$ and
defined by (see Remark \ref{Remark}) a divisor class $C+\alpha$,
$\alpha \in \text{Tor}_d\text{Pic}(X)$, is called {\it
$(d,m)$-canonical} if $C\sim mK_X$ and $d\alpha=0$, {\it pure
$(d,m)$-canonical} if $C\sim mK_X$ and $\alpha =0$, and {\it
numerically $(d,m)$-canonical} if $C\equiv mK_X$. Note that if $f$
is $(d,m)$-canonical, then $B\in |dmK_X|$.

In the following, given a divisor $D$ on a surface $X$ we denote by
$\varphi_{D}:X\to \mathbb P^{\dim |D|}$ the rational map defined by
the complete linear system $|D|$.

The famous Bombieri Theorem (\cite{B}) states that if $m\geq 5$ then
for a smooth minimal (that is, without $(-1)$-curves) projective
surface $X$ the $m$-canonical map $\varphi_{mK_X}: X\to \mathbb
P^{P_m-1}$ is a birational morphism onto its image, where $P_m=\dim
H^0(X, \mathcal O_X(mK_X))$. If to examine the proof of this
theorem, one can see that except Ramanujam vanishing theorem only
the numerical properties of the canonical class $K_X$ are used in
the proof. Therefore Bombiery Theorem is true not only for
$m$-canonical maps, but also for maps $\varphi_{D_m}: X\to \mathbb
P^{P_m-1}$, where $D_m$ is a divisor numerically equivalent to
$mK_X$. In particular, if $dm\geq 5$ and $D\equiv dmK_X$, then by
Bertini Theorem, a generic curve $B\in |D|$ is non-singular and
irreducible.

\begin{prop}\label{2K} Let  $f:Y\to X$ be a numerically $(d,m)$-canonical
totally ramified cyclic covering of a surface $X$. Then $Y$ has the
following invariants:
\begin{equation} \label{pa} \displaystyle p_a(Y)= dp_a(X)+\frac{d(d-1)m((2d-1)m+3)}{12}K^2_X,
\end{equation}

\begin{equation} \label{K2} K^2_Y=d(dm-m+1)^2K^2_X ,\end{equation}
and $q(Y)=q(X)$, where $p_a=p_g-q+1$ is the arithmetic genus of a
surface.
\end{prop}
\proof We have $dK_Y\sim f^*(dK_X+(d-1)B)\equiv f^*(d(dm-m+1)K_X)$.
Therefore $K^2_Y=d(dm-m+1)^2K^2_X$.

It follows from Proposition \ref{pi} that $q(Y)=q(X)$.

Let $e(V)$ denote the topological Euler characteristic of a variety
$V$. We have $$K^2_X+e(X)=12p_a(X)$$ by Noether formula and, by
adjunction formula,
$$-e(B)=2g(B)-2=(B,B+K_X)_X=(dmK_X,(dm+1)K_X)_X=dm(dm+1)K^2_X.$$
Since $f^*(B)=dR$, we have
$$\begin{array}{ll} e(Y)= & d(e(X)-e(B))+e(B)= \\ & d(12p_a(X)-K_X^2+dm(dm+1)K^2_X) -dm(dm+1)K^2_X=
\\ & 12dp_a(X)+d[(d-1)(dm+1)m-1]K^2_X.\end{array}$$
Hence,
$$\begin{array}{ll} e(Y)+K_Y^2= & 12dp_a(X)+d[(d-1)(dm+1)m-1]K^2_X+ d(dm-m+1)^2K^2_X= \\ & 12dp_a(X)+d[(d-1)(dm+1)m-1+
(dm-m+1)^2]K^2_X=
 \\ & 12dp_a(X)+d(d-1)m[(2d-1)m+3]K^2_X
\end{array}$$ and therefore
$$\displaystyle p_a(Y)=dp_a(X)+\frac{d(d-1)m[(2d-1)m+3]}{12}K_X^2 . $$
\qed

By Proposition \ref{2K}, the invariants $K^2_Y$ and $p_a(Y)$ are the
same for all numerically $(d,m)$-canonical totally ramified cyclic
coverings $f:Y\to X$ with fixed $d$ and $m$. We denote them by
$k_{X,d,m}$ and $p_{X,d,m}$, respectively. The moduli space of
surfaces $Z$ with given invariants $K^2_Z =k$ and $p_a(Z)=p$ will be
denoted by $\mathcal M_{k,p}$.

\begin{thm} \label{mod} Let $X$ be a surface of general type.
If there is an element $\alpha\in \text{Tor}_d\text{Pic}(X)$ such
that $\delta(\alpha)$  is not divisible by $d$ in $\text{Tor}(X)$,
then
for each integer $n\geq 
1$ the moduli space $\mathcal M_{k,p}$ with $k=k_{X,d,dn+1}$ and
$p=p_{X,d,dn+1}$ consists of at least two connected components.
\end{thm}
\proof Let us choose a divisor $C\sim (dn+1)K_X$. Consider two
$(d,dn+1)$-canonical totally ramified cyclic coverings $f_i:Y_i\to
X$, $i=1,2$, branched along a smooth irreducible curve $B\in
|d(dn+1)K_X|$, where $f_1$ is the pure $(d,dn+1)$-canonical and
$f_2$ is the $(d,dn+1)$-canonical totally ramified cyclic covering
defined by $C+\alpha$ where $\alpha$ belongs to
$\text{Tor}_d\text{Pic}(X)$ and $\delta(\alpha)$ is not divisible by
$d$ in $\text{Tor}(X)$. For $n\in \mathbb N$ the existence of such a
curve $B$ follows from the inequality $d(dn+1)\geq 6$.

Consider the covering $f_1$. It is given by adding to the field
$\mathbb C(X)$ a function $w_1$ such that $w_1^d =g_1$ and
$(g_1)=B-dC$. Since $$R_1\sim f_1^{*}(C)\sim f^*((dn+1)K_X),$$ then
from the projection formula for the canonical divisors we get
$$K_{Y_1}\sim f_1^*(K_X)+(d-1)R_1\sim f_1^*(K_X)+(d-1)f^*((dn+1)K_X)\sim d(dn-n+1)f^*(K_X).$$
In particular, $K_{Y_1}$ is divisible by $d$ in $\text{Pic}(Y_1)$,
and hence its cohomology class is divisible by $d$ in
$H^2(Y_1,\mathbb Z)$.

The covering $f_2$ is given by adding to the field $\mathbb C(X)$ a
function $w_2$ such that $w_2^d =g_2$ and $(g_2)=B-dC_2$ with
$C_2\sim (dn+1)K_X+\alpha$. Due to Lemma \ref{cov}, $R_2\sim
f_2^{*}(C_2)\sim f_2^*((dn+1)K_X+ \alpha)$. Hence,
$$\begin{array}{ll} K_{Y_2}\sim & f_2^*(K_X)+(d-1)R_2\sim f_2^*(K_X)+(d-1)f_2^*((dn+1)K_X+
\alpha)\sim
\\ & d(dn-n+1)f_2^*(K_X)+(d-1)
f_2^*(\alpha). \end{array}$$ Since $f_2^* :
\text{Tor}(X)\to\text{Tor}(Y_2)$ is an isomorphism and
$\delta(\alpha)\in \text{Tor}(X)$ is not divisible by $d$, this
shows that the cohomology class of $K_{Y_2}$ is not divisible by $d$
in $H^2(Y_2,\mathbb Z)$.

In view of different divisibility of the cohomology classes of their
canonical divisors, $Y_1$ and $Y_2$ can not belong to the same
connected component in the moduli space. \qed

\begin{rem} \label{re} {\rm Theorem \ref{mod} is true also in the case $n=0$ if there is a smooth irreducible
curve $B\in |dK_X|$. From Theorem 5.2 in \cite{B-H-P-V}, Proposition
3 in \cite{Re}, and Bertini Theorem, it follows that a generic curve
$B\in |dK_X|$ is smooth and irreducible if $d\geq 5$, or if $d=4$
and $K_X^2\geq 2$, or if $d=3$ and $K_X^2\geq 3$, or if $d=2$ and
$K_X^2\ge 5$.}
\end{rem}

\begin{rem} \label{rem+} {\rm Note that if
$X_1$ and $X_2$ are two surfaces of general type such that
$K_{X_1}^2=K_{X_2}^2$ and $p_a(X_1)=p_a(X_2)$, but $\pi_1(X_1)$ is
not isomorphic to $\pi_1(X_2)$, and if $f_1:Y_1\to X_1$ and
$f_2:Y_2\to X_2$ are numerically $(d,m)$-canonical totally ramified
cyclic coverings, then, by Proposition \ref{pi}, $Y_1$ and $Y_2$
belong to different connected components of the moduli space
$\mathcal M_{k,p}$ with $k=k_{X_1,d,m}$ and $p=p_{X_1,d,m}$.}
\end{rem}

Let $X$ be a Campedelli surface, that is, $X$ is a surface of
general type with $p_g=0$, $K_X^2=2$, and $\pi_1(X)\simeq (\mathbb
Z/2\mathbb Z)^3$. A totally ramified numerically $(2,1)$-canonical
cyclic covering $f:Y\to X$ is a $(2,1)$-canonical covering since
each non-zero element of  $\text{Tor}(X)\simeq (\mathbb Z/2\mathbb
Z)^3$ is not divisible by two. By Proposition \ref{2K}, a totally
ramified $(2,1)$-canonical cyclic covering $Y$ of $X$ has the
following invariants: $K_Y^2=16$ and $p_a=4$. Note also that, by
Lemma \ref{cov}, for given non-singular curve $B\in |2K_X|$ there
exist eight different totally ramified $(2,1)$-canonical cyclic
coverings of $X$ branched along $B$.
\begin{prop} \label{Camp1} There are exactly two connected components of the moduli
space $\mathcal M_{16,4}$ to which belong the totally ramified
$(2,1)$-canonical cyclic coverings of Campe\-delli surfaces.
\end{prop}
\proof Recall that every Campedelli surface $X$ can be obtained as
follows (see, for example, \cite{K}). Let $\widetilde L=L_1\cup
\dots \cup L_7$ be a line arrangement in $\mathbb P^2$ consisting of
seven lines. We numerate them by the non-zero elements
$\alpha_i=(a_{i,1},a_{i,2},a_{i,3})\in (\mathbb Z/2\mathbb Z)^3$ and
make two assumptions: the arrangement $\widetilde L$ has no $r$-fold
points with $r\geq 4$; and if $\widetilde L$ has a tripe point
$p_{\alpha_{i_1},\alpha_{i_2},\alpha_{i_3}}= L_{\alpha_{i_1}}\cap
L_{\alpha_{i_2}}\cap L_{\alpha_{i_3}}$, then
$\alpha_{i_1}+\alpha_{i_2}+\alpha_{i_3}\neq 0$.  Consider the Galois
covering $\widetilde g: \widetilde X\to \mathbb P^2$ with Galois
group $\text{Gal}(\widetilde X/\mathbb P^2)\simeq (\mathbb
Z/2\mathbb Z)^3$, that is, branched in $\widetilde L$ and defined by
the epimorphism $\varphi :H_1(\mathbb P^2\setminus \widetilde L,
\mathbb Z)\to G=(\mathbb Z/2\mathbb Z)^3$ given by $\varphi
(\lambda_{i})=\alpha_i$, where $\lambda_{i}$ is the element of
$H_1(\mathbb P^2\setminus \widetilde L, \mathbb Z)$ represented by a
small circle around $L_{\alpha_i}$.

The only singular points of $\widetilde X$ are the points lying over
the triple points $p_{\alpha_{i_1},\alpha_{i_2},\alpha_{i_3}}$, and
$X$ that is the resolution of singularities of $\widetilde X$ is a
Campedelli surface. To resolve the singularities of $\widetilde X$,
it suffices to blow up the triple points of $\widetilde L$.  The
composition $\sigma : \widetilde{\mathbb P}^2\to {\mathbb P}^2$ of
blowups with centers at the triple points of $\widetilde L$ induces
the Galois covering $g: X \to \widetilde{\mathbb P}^2$ and  we have
(\cite{K}) $2K_X \sim g^*(L)$, where $L$ is a line in $\mathbb P^2$.
(If $\widetilde L$ has no triple points, then $g=\widetilde g$.)

Let $X$ be defined by a Campedelli arrangement $\widetilde L$ and
let $f$ be branched along $B=g^{-1}(L_8)$, where a line
$L_8\not\subset\widetilde L$. Without loss of generality, we can
assume that the line arrangement $\mathcal L=\widetilde L\cup L_8$
is generic. Then the fundamental group $\pi_1(\mathbb P^2\setminus
\mathcal L)$ is abelian and hence $h=g\circ f:Y\to\mathbb P^2$ is a
$(\mathbb Z/2\mathbb Z)^4$-Galois covering of $\mathbb P^2$ branched
along $\mathcal L$ and defined by an epimorphism $\psi :H_1(\mathbb
P^2\setminus \mathcal L, \mathbb Z)\to (\mathbb Z/2\mathbb Z)^4$
given by $\psi(\lambda_{i})=
\widetilde{\alpha}_i=(a_{i,1},a_{i,2},a_{i,3},a_{i,4})$ for some
$a_{i,4}\in \mathbb Z/2\mathbb Z$, $i=1,\dots,7$, and
$\psi(\lambda_8)=\widetilde{\alpha}_8=(0,0,0,1)$, where $\lambda_8$
is the element  represented by a small circle around $L_8$. Note
that $\psi$ is defined by $h$ up to automorphisms of $(\mathbb
Z/2\mathbb Z)^4$ (that is, by the choice of the basis in $(\mathbb
Z/2\mathbb Z)^4$). We have $\sum_{i=1}^8\lambda _i=0$ in
$H_1(\mathbb P^2\setminus \mathcal L, \mathbb Z)$. Therefore the
subset $M=\{ \widetilde{\alpha}_i \subset \mathbb P^3_{\mathbb
Z/2\mathbb Z}\mid i=1,\dots,8\}$ of the projective space $\mathbb
P^3=\mathbb P^3_{\mathbb Z/2\mathbb Z}$ over the field $\mathbb
Z/2\mathbb Z$ is  {\it totally even,} that is, it satisfies the
following property:
{\it for each plane $\mathbb P^2\subset\mathbb P^3$ the intersection
$M\cap \mathbb P^2$ consists of even number of points.}

\begin{lem} \label{ev} Up to the action of $\text{PGL}(4,\mathbb Z/2\mathbb
Z)$, there exist exactly two different totally even subsets $M$ of
the projective space $\mathbb P^3$ over the field $\mathbb
Z/2\mathbb Z$ consisting of eight points, namely, either $M=\mathbb
A^3$ (type $I$), or
$$\begin{array}{ll} M=\{ & (1,0,0,0), (1,1,0,0), (0,1,0,0), (1,0,1,0), (0,0,1,0), \\ & (0,1,1,0), (1,1,1,1), (0,0,0,1)\, \, \} \end{array}$$
(type $II$). \end{lem} \proof Type $I$ corresponds to the case when
there is a plane $P\subset \mathbb P^3$ (say, $P$ is given by
equation $a_4=0$) such that $M\cap P=\emptyset$ and in this case
$$M=\mathbb A^3= \{ \widetilde{\alpha}=(a_1,a_2,a_3,a_4)\in \mathbb P^3\mid a_4=1 \} .$$

Assume that there is not a plane $P$ such that $M\cap P=\emptyset$.
For each plane $P_i\subset \mathbb P^3$ denote by $n_i=|M\cap P_i|$
the number of points of the intersection $M\cap P_i$. Each number
$n_i$ is even, $0<n_i\leq 6$.

First of all, let us show that the condition $n_i=6$ for some $i$
(say, $n_1=6$) is equivalent to the condition that there exists a
line  $l\subset M$ (and in this case $|M\cap l|=3$). Indeed, if
$n_1=6$ then there is a unique point $\widetilde{\alpha}_0\in P_1$
such that $\widetilde{\alpha}_0\not\in M$, since $|P_1|=7$.
Therefore for any line $l\subset P_1$ not passing through
$\widetilde{\alpha}_0$ we have $l\subset M$. Conversely, for a line
$l\subset M$ consider the pencil of planes passing through $l$. It
consists of three planes, say, $P_1$, $P_2$, and $P_3$. We have
$n_i\geq 4$ for $i=1,2,3$ and $(n_1-3)+(n_2-3)+(n_3-3)=8-3=5$, that
is, $n_1+n_2+n_3=14$. Up to permutation, this equation has the only
one solution consisting of even numbers
and satisfying inequalities $n_i\geq 4$, namely, $n_1=6$ and
$n_2=n_3=4$.

Let us show that $M$ is of type $II$ if there is a plane $P$ such
that $|M\cap P|=6$. Without loss of generality, we can assume that
$P$ is given by equation $a_4=0$ and the points
$\widetilde{\alpha}_7$ and $\widetilde{\alpha}_8$ lying in
$M\setminus P$ have coordinates: $\widetilde{\alpha}_7=(1,1,1,1)$
and $\widetilde{\alpha}_8=(0,0,0,1)$. Since $|M\cap P|=6$ and
$|P|=7$, to prove that $M$ is of type $II$, it suffices to show that
the point $\widetilde{\alpha}_0=(1,1,1,0)\not\in M$. Assume that
$\widetilde{\alpha}_0=(1,1,1,0)\in M$. Then there is another point
$\widetilde{\alpha}_1=(a_{1,1}, a_{1,2},a_{1,3},0)\not\in M$ with
$a_{1,i}=0$ for some $i=1,2,3$. Without loss of generality, applying
a projective transformation permuting coordinates and preserving $P$
fixed, we can assume that this point is $\widetilde{\alpha}_1=(0,
a_{1,2},a_{1,3},0)$. Then $|M\cap l|=2$, where the line $l\subset
\mathbb P^3$ is given by equations $a_1=a_4=0$. Therefore we have
the only three points belonging to $M$ whose first coordinate $a_1$
vanishes, namely, two points lying in $M\cap l$ and the point
$(0,0,0,1)$. Therefore $|M\cap P_1|=3$, where the plane
$P_1=\{\widetilde{\alpha}=(a_1,a_2,a_3,a_4)\in \mathbb P^3\mid a_1=0
\}$. Contradiction.

Let us show that the case $0<n_i\leq 4$ for all planes $P_i\subset
\mathbb P^3$ is impossible. First of all, assume that $n_i=4$ for
all planes $P_i$. On the one hand, the number of planes in $\mathbb
P^3$ is equal to $15$. On the other hand, the number of planes
passing through a point $\widetilde{\alpha}\in \mathbb P^3$ is equal
to $7$ and the inequality $4\cdot 15\neq 8\cdot 7$ implies that this
case is impossible.

Finally, assume that $0<n_i\leq 4$ for all planes $P_i\subset
\mathbb P^3$ and there is a plane $P_1$ such that $n_1=|M\cap
P_1|=2$. Then there exists a line $l\subset P_1$ such that $M\cap
l=\emptyset$. Consider the pencil of planes passing through $l$. It
consists of three planes, $P_1$, $P_2$, and $P_3$. Since
$n_1+n_2+n_3=8$, we have $n_2+n_3=6$. Therefore, we can assume that
$n_2=2$ and $n_3=4$. We have $M\cap P_3=P_3\setminus l$ and hence
each plane $P$ not containing $l$ has two common points with $P_3$
belonging to $M$. On the other hand, let
$\{\widetilde{\alpha}_1,\widetilde{\alpha}_2\}=M\cap P_1$ and
$\{\widetilde{\alpha}_3,\widetilde{\alpha}_4\}=M\cap P_2$. Then the
plane $P_4$ passing through the points $\widetilde{\alpha}_1$,
$\widetilde{\alpha}_2$, and $\widetilde{\alpha}_3$ does not contain
the line $l$. Therefore $n_4\geq 5$. Contradiction.
\qed \\

Let us return to the proof of Proposition \ref{Camp1}. Let the lines
of a line arrangement $\mathcal L$ consisting of eight lines be
numerated by the points of a totally even set $M\subset (\mathbb
Z/2\mathbb Z)^4$, $\mathcal L=\mathcal
L_M=\bigcup_{\widetilde{\alpha}_i\in M}L_{\widetilde{\alpha}_i}$.
The line arrangement $\mathcal L_M$ defines an epimorphism
$\psi:H_1(\mathbb P^2\setminus \mathcal L_M, \mathbb Z)\to (\mathbb
Z/2\mathbb Z)^4$ given by $\psi(\lambda_{\widetilde{\alpha}_i})=
\widetilde{\alpha}_i$ and the epimorphism $\psi$ defines a $(\mathbb
Z/2\mathbb Z)^4$-Galois covering $h:Y\to \mathbb P^2$ branched along
$\mathcal L_M$. Since the epimorphism $\psi$ is defined by $h$ only
up to automorphism of $(\mathbb Z/2\mathbb Z)^4$, therefore by Lemma
\ref{ev}, we can assume that $M$ is either $M=\{
\widetilde{\alpha}=(a_1,a_2,a_3,a_4)\in \mathbb P^3\mid a_4=1 \} $
(type $I$) or
$$\begin{array}{ll} M=\{ & (1,0,0,0), (1,1,0,0), (0,1,0,0), (1,0,1,0), (0,0,1,0), \\ & (0,1,1,0), (1,1,1,1), (0,0,0,1)\, \, \} \end{array}$$
(type $II$). Now, to complete the proof, by Remark \ref{re}, it
suffices to note that the set of line arrangements $\mathcal L_M$,
where $M$ is of type $I$ (respectively, of type $II$), is an
everywhere dense Zariski open (and, consequently, connected) subset
of $(\mathbb P^2)^8$. \qed

\begin{thm}\label{def-eq} Let $X$ be a surface of general type, and let $f_1:Y_1\to
X$, $f_2:Y_2\to X$ be two numerically $(d,m)$-canonical totally
ramified cyclic coverings defined, respectively, by divisors
$C_1\sim mK_X+ \alpha_1$ and $C_2\sim mK_X+ \alpha_2$, where $
\alpha_1$ and $\alpha_2$ are numerically equivalent to zero and such
that $\delta(\alpha_1)=\delta(\alpha_2)$. If $dm\geq 5$ then $Y_1$
and $Y_2$ are deformation equivalent.
\end{thm}
\proof Let $B_i\in |dC_i|$ be the branch curve of the covering
$f_i$, $i=1,2$. The covering $f_i$ is given by adding a function
$w_i$ to the field $\mathbb C(X)$ with $w_i^d =g_i$,
$(g_i)=B_i-dC_i$, $C_i\sim mK_X+ \alpha_i$. Let
$\Delta=\delta(\alpha_1)=\delta(\alpha_2)\in\text{Tor}(X)$. Denote
by $\text{Pic}_\Delta(X)=\delta^{-1}(\Delta)\subset \text{Pic}(X)$.

Consider the scheme $T_{\Delta, dm}$ parametrizing the curves $B_t$
in $X$ numerically equivalent to $dmK_X$ and such that $\delta(B_t)=
dm\delta(K_X)+d\Delta$ (the scheme $T_{\Delta,dm}$ is fibered over
$\text{Pic}_{\Delta}(X)$, $\gamma_{1}:
T_{\Delta,dm}\to\text{Pic}_{\Delta}(X)$, with fibers $\gamma_1^{-1}(
\alpha_t)=\mathbb P(H^0(X, \mathcal O_X(dmK_X+ d \alpha_t))$ over $
\alpha_t\in \text{Pic}_{\Delta}(X)$; see \cite{Mu}). Obviously, the
subscheme $U$ consisting of points $t$ of $T_{\Delta,dm}$ for which
$B_t$ are smooth curves is a Zariski open non-empty subset.

Similarly, let $T_{\Delta, 2}$ be the scheme parametrizing the
curves $D_t$ in $X$ numerically equivalent to $2K_X$ and such that
$\delta(D_t)= 2\delta(K_X)+\Delta$. The scheme $T_{\Delta,2}$ also
is fibered over $\text{Pic}_{\Delta}(X)$, $\gamma_{2}:
T_{\Delta,2}\to\text{Pic}_{\Delta}(X)$, with fibers
$\gamma_2^{-1}(\alpha_t)=\mathbb P(H^0(X, \mathcal O_X(2K_X+
\alpha_t))$ over $
\alpha_t\in \text{Pic}_{\Delta}(X)$. Denote by
$T=U\times_{\text{Pic}_{\Delta}(X)}T_{\Delta,2}$ the product of the
fibrations $\gamma_1$ and $\gamma_2$ and let $p_i$, $i=1,2$, be the
projections of $T$ onto the factors.

Let us fix a divisor $D\sim (m-2)K_X$ and associate with each $t\in
T$ the divisor $C_{p_2(t)}=D_{p_2(t)}+D\sim mK_X+\alpha_{p_2(t)}$.
By Lemma \ref{cov}, the coverings $f_1$ and $f_2$ are defined by
divisors $C_{p_2(t_i)}=D_{p_2(t_i)}+D\sim mK_X+\alpha_{p_2(t_i)}$
for some points $t_i\in T$, $i=1,2$, and branched along the curves
$B_{p_1(t_i)}$. The family of divisors $D_t=B_{p_1(t)}-dC_{p_2(t)}$
defines a divisor $\widetilde D$ in $X\times T$ such that
$\widetilde D\cap (X\times \{ t\})= D_t$ for each $t\in T$. By
Corollary 6 in \cite{Mu1}, there is an invertible sheaf $\mathcal L$
on $T$ such that $\mathcal O_{X\times T}(\widetilde D)=p^*(\mathcal
L)$, where $p:X\times T\to T$ is the projection. Obviously, we can
choose a rational section of $\mathcal L$ the support of whose
divisor $L$ does not contain the points $t_1$ and $t_2$. Let us
consider a rational function $\widetilde g=g_t$, $t\in T$, $t\not\in
\text{Supp}\, L$, such that the divisor $(\widetilde g)$ is
$\widetilde D-p^*(L)$. The function $\widetilde g$ defines a cyclic
covering $\widetilde f:\widetilde Y\to X\times T$ given by
$\widetilde w^d=\widetilde g$. This covering can be considered as a
connected family of cyclic coverings $f_t:Y_t\to X$, $t\in
T\setminus \text{Supp}\, L$ given by $C_{p_2(t)}=D_{p_2(t)}+D\sim
mK_X+\alpha_{p_2(t)}$ and  branched along the curves $B_{p_1(t_i)}$,
and hence $Y_1$ and $Y_2$ are deformation equivalent. \qed

\begin{rem} \label{rem}{\rm Theorem \ref{cor-main2} is true without assumption $dm\geq 5$ if
for each $\alpha\in \text{Pic}_{\Delta}$ a generic curve $B\sim
d(mK_X+\alpha)$ is smooth.}\end{rem}

\begin{thm} \label{pg0}
Let $X$ be a surface of general type with $p_g=0$, and let $f:Y\to
X$ be a numerically $(2,m)$-canonical totally ramified cyclic
covering of $X$. Then the rational map $\varphi_{K_Y}: Y\to \mathbb
P^{p_g(Y)-1}$  factorizes through $f$, that is, there exists a
rational map $\psi: X\to \mathbb P^{p_g(Y)-1}$ such that
$\varphi_{K_Y}=\psi\circ f$.
\end{thm}
\proof Let $B \equiv 2mK_X $ be the branch curve of $f$ and let $f$
be defined by a divisor $C\equiv mK_X$. As in the proof of Theorem
\ref{mod}, due to the projection formula for the canonical divisor
we have the equality
\begin{equation} \label{f1} K_Y=f^*(K_X+C).\end{equation}
As any surface of general type with $p_g=0$, the surface $X$ is
regular (that is, $q(X)=0$), so that Mumford Vanishing Theorem and
Riemann-Roch Theorem imply
\begin{equation} \label{f2} \dim H^0(X, \mathcal
O_V(K_X+C))=\frac{(K_X+C,C)_X}{2}+1=\frac{m(m+1)}{2}K^2_X+1.\end{equation}
By Proposition \ref{pi}, the surface $Y$ is also regular; therefore,
applying the formula (\ref{pa}) from Proposition \ref{2K} we obtain
\begin{equation} \label{f3}
\displaystyle p_g(Y)=1+\frac{m(m+1)}{2}K_X^2 .\end{equation} Now,
the theorem follows from (\ref{f1}) -- (\ref{f3}). \qed

\begin{cor} \label{pg01}
Let $X$ be a surface of general type with $p_g=0$, and  let $f:Y\to
X$ be a numerically $(2,m)$-canonical totally ramified cyclic
covering of $X$ defined by a divisor $C\equiv mK_X$. Then for $m\geq
4$ the canonical map $\varphi_{K_Y}: Y\to \mathbb P^{p_g(Y)-1}$ is a
morphism of degree $2$ over its image, and this image,
$\varphi_{K_Y}(Y) \subset \mathbb P^{p_g(Y)-1}$, coincides with the
image of $X$ under the birational morphism $\psi=\varphi_{K_X+C}$.
\qed
\end{cor}

In the present time, there is a long list of known surfaces of
general type with $p_g=0$ (see, for example, the survey
\cite{B-C-P}), but up to now, this list is far from completeness. In
the following three propositions we investigate the degree of the
canonical map of a pure $(2,1)$-canonical totally ramified cyclic
covering, respectively, of a Campedelli surface, a Burniat surface,
and a Mendes Lopes -- Pardini surface.

\begin{prop} \label{Camp}  Let $f:Y\to
X$ be a pure $(2,1)$-canonical totally ramified cyclic covering of a
Campedelli surface $X$. Then $\varphi_{K_Y}:Y\to\mathbb P^2$ is a
morphism of degree $\deg \varphi_{K_Y}=16$.
\end{prop}
\proof In notation used in the proof of Theorem \ref{Camp1}, we have
(\cite{K}) $\varphi_{2K_X}=\sigma\circ g$. Therefore, $\deg
\varphi_{2K_X}=8$. Now, Theorem \ref{pg0} implies that
$\varphi_{K_Y}=\varphi_{2K_X}\circ f$, and thus the map
$\varphi_{K_Y}$ is regular and of degree $\deg \varphi_{K_Y}=16$.
\qed

\begin{prop} \label{Burn}  Let $f:Y\to X$ be a
pure $(2,1)$-canonical totally ramified cyclic covering of a Burniat
surface $X$ with $K_X^2\geq 3$. Then $\varphi_{K_Y}:Y\to \mathbb
P^{K_X^2}$ is a morphism of degree $8$ over its image, and this
image, $Z=\varphi_{K_Y}(Y)$, is a Del Pezzo surface with
$K^2_Z=K_X^2$ embedded into $\mathbb P^{K_X^2}$ by means of the
anticanonical map.
\end{prop}
\proof Similar to Campedelli surfaces, a Burniat surface $X$ is the
resolution of singularities  of a Galois covering $\widetilde g:
\widetilde X\to \mathbb P^2$ with Galois group
$\text{Gal}(\widetilde X/\mathbb P^2)\simeq (\mathbb Z/2\mathbb
Z)^2$ branched along a Burniat line arrangement $\widetilde L$ (see
details in \cite{Bu} or in \cite{K}). To resolve the singularities,
it suffices to blow up the $r$-fold points of the line arrangement
$\widetilde L$ with $r\geq 3$ and consider the induced Galois
covering $g: X \to \widetilde{\mathbb P}^2$, where $\sigma :
\widetilde{\mathbb P}^2\to {\mathbb P}^2$ is the composition of
blowups with centers at the $r$-fold points of $\widetilde L$,
$r\geq 3$. The number of $r$-fold points of $\widetilde L$ with
$r\geq 3$ is equal to $9-K^2_X$ and therefore $\widetilde{\mathbb
P}^2$ is a Del Pezzo surface with $K^2_{\widetilde{\mathbb
P}^2}=K_X^2$. Moreover, $2K_X=g^*(-K_{\widetilde{\mathbb P}^2})$ and
$\dim H^0(X,\mathcal O_X(2K_X))=\dim H^0(\widetilde{\mathbb
P}^2,\mathcal O_{\widetilde{\mathbb P}^2}(-K_{\widetilde{\mathbb
P}^2}))$ (see, for example, \cite{K}). Therefore
$\varphi_{2K_X}=\varphi_{-K_{\widetilde{\mathbb P}^2}}\circ g$ is a
regular map of degree four. Now, Theorem \ref{pg0} shows that
$\varphi_{K_Y}=\varphi_{2K_X}\circ f$,
which implies that $\deg \varphi_{K_Y}=8$. \qed \\

Mendes Lopes and Pardini (\cite{M-P}) constructed a six-dimensional
family of surfaces $X$ of general type with $p_g=0$ and $K^2_X=3$.
For each of these surfaces, call them {\it Mendes Lopes -- Pardini
surfaces}, the map $\varphi_{2K_X}$ is  regular of degree $2$ over
its image, $Z$, which is a singular Enriques surface $Z\subset
\mathbb P^3$, $\deg Z=6$. We apply Theorem \ref{pg0} following the
same lines as above and obtain the following result.

\begin{prop} \label{Me-Pa}  Let $f:Y\to X$ be a
pure $(2,1)$-canonical totally ramified cyclic covering of a Mendes
Lopes -- Pardini surface $X$. Then $\varphi_{K_Y}:Y\to \mathbb
P^{3}$ is  a regular map of degree $4$ over its image, $Z$, which is
a singular Enriques surface $Z\subset \mathbb P^3$, $\deg Z=6$. \qed
\end{prop}

\section{Cyclic coverings of rigid surfaces}
Recall that each Miyaka-Yau surface $X$ being a holomorphic quotient
of a ball satisfies Mostow rigidity theorem. The latter, in one of
its well known formulations, states that (existence part) each
element of the group $\operatorname{Out} \pi_1(X)$ is realized by an
antiholomorphic or holomorphic diffeomorphism $g:X\to X$, and
furthermore (unicity part) such a realization is unique ({\it cf.,}
\cite{Mo} and \cite{GP}).

Given a complex surface $X$, denote by $\bar X$ the surface with the
complex conjugate complex structure and by $\text{Kl} (X)$ the group
of holomorphic and anti-holomorphic automorphisms of $X$.

\begin{thm} \label{main2} Let $X$ be a
Miyaoka-Yau surface and let $f_1:Y_1\to X$, $f_2:Y_2\to X$ be two
numerically $(d,m)$-canonical totally ramified cyclic coverings
defined, respectively, by divisors $C_1\sim mK_X+ \alpha_1$ and
$C_2\sim mK_X+ \alpha_2$, where both $ \alpha_1$, $\alpha_2$ are
numerically equivalent to zero. Assume, in addition, that in the
case $d\geq 3$ both $\delta(\alpha_1)$  and $\delta(\alpha_2)$ have
the same order $n$ in $\text{Tor}(X)$ coprime with $d-1$. If $Y_1$
and $Y_2$ are orientation preserving diffeomorphic, then there
exists a {\rm (}holomorphic or anti-holomorphic{\rm )} automorphism
$\psi\in \text{Kl}(X)$ such that
$\psi^*(\delta(\alpha_2))=\pm\delta(\alpha_1)$ {\rm (}the sign is
plus if $\psi$ is holomorphic, and minus otherwise{\rm )}.

If $Y_1$ and $Y_2$ are deformation equivalent, then there exists an
automoprhism $\psi\in\text{Aut} (X)$ such that
$\psi^*(\delta(\alpha_2))=\delta(\alpha_1)$.
\end{thm}
\proof Let $\varphi:Y_1\to Y_2$ be an orientation preserving
diffeomorphism.

We have
$$\begin{array}{ll} K_{Y_1}= & f_1^*(K_X)+(d-1)R_1\sim
f_1^*(K_X)+(d-1)f_1^*(mK_X+ \alpha_1)\sim \\ & (dm-d+1)f_1^*(K_X)
+(d-1) f_1^*(\alpha_1)\end{array}$$ and, similarly, $K_{Y_2}
\sim(dm-d+1)f_2^*(K_X) +(d-1) f_2^*(\alpha_2)$. The surfaces $Y_1$
and $Y_2$ being coverings of a surface of Kodaira dimension $2$ have
also Kodaira dimension $2$, and thus they are of general type. They
are both minimal, since for any curve $E\subset Y_i$ we have $E\cdot
K_{Y_i}=f_*(E)\cdot (dm-d+1)K_X>0$. Therefore, $\pm \delta(K_{Y_i})$
are the only Seiberg-Witten basic classes in $H^2(Y_i;\mathbb Z)$
(see \cite{W} for the case $p_g>0$ and, for example, \cite{Br} for
$p_g=0$), which implies $\varphi^*(\delta(K_{Y_2}))=\pm
\delta(K_{Y_1})$.

By Proposition \ref{pi} the map $\varphi$ induces an automorphism
$\varphi_* :\pi_1(X)=\pi_1(Y_1)\to\pi_1(Y_2)=\pi_1(X)$ of $\pi_1(X)$
(more precisely, an element of $\operatorname{Out}\pi_1(X)$), and
therefore, in accordance with Mostow rigidity, it is induced by an
either holomorphic or antiholomorphic automorphism $\psi:X\to X$. By
Siu rigidity (see \cite{Siu}), the map $f_2\circ \varphi:Y_1\to X$
is homotopic to either holomorphic or antiholomorphic morphism
$\widetilde f_1:Y_1\to X$ and, as it follows once more from Mostow
rigidity, $\widetilde f_1=\psi\circ f_1$ since these two morphisms
define the same element of $\operatorname{Out}\pi_1(X)$.

If $\psi$ is holomorphic then $\psi^*(K_X)=K_X$. Therefore
$\varphi^*(f_2^*(\delta(K_X)))=f^*_1(\delta(K_X))$ and hence
$\varphi^*((d-1)\delta(\alpha_2))=(d-1)\delta(\alpha_1)$. Let us
show that then $\varphi^*(\delta(\alpha_2))=\delta(\alpha_1)$.
Indeed, assume that
$\varphi^*(\delta(\alpha_2))=\delta(\alpha_1)+\beta$ for some
$\beta\neq 0$ and such that $(d-1)\beta=0$. A priori, it is possible
only if $d\geq 3$. But, in this case, since $(n,d-1)=1$, the element
$\delta(\alpha_1)+\beta$ must have the order greater than $n$. On
the other hand, the element $\varphi^*(\delta(\alpha_2))$ has the
same order as $\delta(\alpha_2)$. Therefore $\beta=0$ and
$\psi^*(\delta(\alpha_2))=\delta(\alpha_1)$.

If $\psi$ is antiholomorphic then $\psi^*(K_X)=-K_X$. Thus
$\varphi^*(f_2^*(\delta(K_X)))=-f^*_1(\delta(K_X))$. Hence
$\varphi^*((d-1)\delta(\alpha_2))=-(d-1)\delta(\alpha_1)$ and as
above we obtain that $\psi^*(\delta(\alpha_2))=-\delta(\alpha_1)$.

If $Y_1$ and $Y_2$ are deformation equivalent, then there is an
orientation preserving diffeomorphism $\phi : Y_1\to Y_2$ such that
$\varphi^*(\delta(K_{Y_2}))= \delta(K_{Y_1})$. Repeating once more
the same arguments as above, we conclude that there exists
$\psi\in\text{Aut}(X)$ such that
$\psi^*(\delta(\alpha_2))=\delta(\alpha_1)$.
\qed \\

\begin{cor}\label{cor-main2} Let $f_1:Y_1\to X$ and
$f_2:Y_2\to X$ be numerically $(d,m)$-canonical totally ramified
cyclic coverings as in Theorem {\rm \ref{main2}}. Suppose that
$dm\geq 5$ and there exists an orientation preserving diffeomorphism
$\varphi: Y_1\to Y_2$. Then $Y_2$ is deformation equivalent to $Y_1$
if $\varphi^*(K_{Y_2})=K_{Y_1}$, and $\overline Y_2$ is deformation
equivalent to $Y_1$ if $\varphi^*(K_{Y_2})=-K_{Y_1}$.
\end{cor}
\proof Let $B_i\in |dC_i|$ be the branch curve of the covering
$f_i$, $i=1,2$. The covering $f_i$ is given by adding a function
$w_i$ to the field $\mathbb C(X)$ with $w_i^d =g_i$,
$(g_i)=B_i-dC_i$, $C_i\sim mK_X+ \alpha_i$.

According to Theorem \ref{main2}, there exists $\psi\in\text{Kl}(X)$
such that $\psi^*(\alpha_2)=\pm\alpha_1$ with sign "plus" if $\psi$
is holomorphic and sing "minus" otherwise. If $\psi$ is
antiholomorphic, define an automorphism $\psi^!:
\text{Div}(X)\to\text{Div}(X )$ by $\psi^!(D)=\psi^{-1}(D)$ in the
case $D\subset X$ is a curve.

Let $\psi$ be a holomorphic automorphism. Then the covering
$\psi^{-1}\circ f_2:Y_2\to X$ is given by adding a function
$\widetilde w_2$ to the field $\mathbb C(X)$ with $\widetilde w_2^d
=\psi^*(g_2)$, $(\psi^*(g_2))=\psi^*(B_2)-d\psi^*(C_2)$, where
$\psi^*(C_2)\sim mK_X+\psi^*(\alpha_2)$. By Theorem \ref{def-eq},
$Y_1$ and $Y_2$ are deformation equivalent, since
$\delta(\psi^*(\alpha_2))= \delta(\alpha_1)$.

Let $\psi$ be an anti-holomorphic automorphism. Then the covering
$\psi^{-1}\circ f_2:\overline Y_2\to X$ is given by adding a
function $\overline w_2$ to the field $\mathbb C(X)$ with $\overline
w_2^d =\psi^*(\overline g_2)$,
$(\psi^!(g_2))=\psi^!(B)-d\psi^!(C_2)$, where $\psi^!(C_2)\sim
mK_X+\psi^!(\alpha_2)$. By Theorem \ref{def-eq}, $Y_1$ and
$\overline Y_2$ are deformation equivalent, since
$\delta(\psi^!(\alpha_2))= \delta(-\psi^*(\alpha_2))=
\delta(\alpha_1)$. \qed

\begin{cor}\label{Cnew}
Let $X$ be a Miyaoka-Yau surface and let $d\ge 2$, $m\ge 1$ be
integers such that $dm\ge 5$ and $d-1$ is prime with respect to the
order of the group $\text{Tor}(X)$. Then in the moduli space of
surfaces the number of connected components that contain numerically
$(d,m)$-canonical totally ramified cyclic coverings of $X$ is equal
to the number of orbits of the action of $\text{Aut}(X)$ on
$\text{Tor}(X)$.
\end{cor}

\proof Let $f_1 :Y_1\to X_1$, $f_2 :Y_2\to X_2$ be two totally
ramified numerically $(d,m)$-canonical cyclic coverings given by
divisors $C_1\sim mK_X+\alpha_1$, $C_2\sim mK_X+\alpha_2$,
respectively, where $\alpha_1$ and $\alpha_2$ are both numerically
equivalent to zero. Assume that there exists an automorphism
$\psi\in\text{Aut}(X)$ such that
$\psi^*(\delta(\alpha_2))=\delta(\alpha_1)$. Then, $\psi^{-1}\circ
f_2 : Y_2\to X$ is a totally ramified numerically $(d,m)$-canonical
cyclic covering given by divisor $\psi^*(C_2)\sim
mK_X+\psi^*(\alpha_2)$ with
$\delta(\psi^*(\alpha_2))=\psi^*(\delta(\alpha_2))=\delta(\alpha_1)$,
and therefore, according to Theorem \ref{def-eq}, the surfaces $Y_1$
and $Y_2$ belong to the same connected component of the moduli
space.

The reverse statement follows from Theorem \ref{main2}. \qed

\begin{thm}\label{newTh} Let $X$ be a Miyaoka-Yau surface and let $f_1:Y_1\to X$, $f_2:Y_2\to \bar X$ be
two totally ramified numerically $(d,m)$-canonical cyclic coverings.
If $Y_1$ and $Y_2$ are deformation equivalent, then $\text{Kl}(X)\ne
\text{Aut}(X)$.
\end{thm}

\proof The covering $f_2$ provides also a totally ramified
numerically $(d,m)$-canonical cyclic covering $f_2:\bar Y_2\to X$.
As in the proof of Theorem \ref{main2} we have
$$
K_{Y_1}\equiv (dm-d+1)f^*_1(K_X) \quad\text{and}\quad K_{\bar
Y_2}\equiv (dm-d+1)f^*_2(K_X).
$$

Since $Y_1$ and $Y_2$ are deformation equivalent, there exists an
orientation preserving diffeomorphism $\phi:Y_1\to Y_2$ such that
$\phi^*(\delta(K_{\bar Y_2}))=-\delta(K_{Y_1})$. Therefore, we have
\begin{equation} \label{equ} \varphi^*(f_2^*(\delta(K_{X})))\equiv -f_1^*(\delta(K_X)).\end{equation}

The same arguments as in the proof of Theorem \ref{main2} shows that
the group automorphism $\phi_*:
\pi_1(X)=\pi_1(Y_1)\to\pi_1(Y_2)=\pi_1(X)$ (more precisely, the
element of the group $\text{Out}\pi_1(X)$) is induced by a
holomorphic, or anti-holomorphic, map $\psi : X\to X$ such that the
map $f_2\circ\phi : Y_1\to X$ is homotopic to $\psi\circ f_1 :
Y_1\to X$, since the both maps define the same element of
$\text{Out}\pi_1(X)$.

Let us show that $\psi$ is
an anti-holomorphic automorphism. Indeed, if $\psi$ is holomorphic,
then $\psi^*(\delta(K_X))=\delta(K_X)$. Therefore, it follows from
equality (\ref{equ}) that
$$-f_1^*(\delta(K_X))\equiv
\varphi^*(f_2^*(\delta(K_X)))=f^*_1(\psi^*(\delta(K_X)))=f^*_1(\delta(K_X)),$$
but it is impossible since the element $f^*_1(\delta(K_X))$ is not
numerically equivalent to zero in $H^2(Y_1,\mathbb Z)$.\qed

\begin{cor} \label{Cplus} For each pair of positive integers $d, m$ with $dm\ge
5$, $d\not \equiv 1 (\text{mod}\, 5)$ and a surface $X$ of  general
type with $(K^2)_{X}=333$ and $e(X)=111$, the moduli space $\mathcal
M_{k_{X,d,m},p_{X,d,m}}$ has at least $3\cdot 5^6$ different
connected components.
\end{cor}
\proof Two surfaces $X$ as in the statement are constructed in
Examples 1 and 2 in \cite{KK1}, denote them by $X_1$ and $X_2$
respectively. They are both obtained by resolution of singularities
of the abelian $(\mathbb Z/5\mathbb Z)^2$-coverings $\widetilde
h_i:\widetilde X_i\to\mathbb P^2$ branched along the line
arrangement $\widetilde L=\cup_{j=1}^9L_j$ dual to the inflection
points of a smooth plane cubic. To resolve the singularities of
$\widetilde X_i$, we blow up the $3$-fold  points of $\widetilde L$
and take the normal closure of $\widetilde{\mathbb P}^2$ in the
field $\mathbb C(\widetilde X_i)$; denote this blow up by
$\sigma:\widetilde{\mathbb P}^2\to \mathbb P^2$ and the induced
covering by $h_i:X_i\to\widetilde{\mathbb P}^2$.

The surfaces $\widetilde X_i$ can be obtained also as factor-spaces
$Z/\widetilde G_i$, where $\widetilde G_i\simeq (\mathbb Z/5\mathbb
Z)^6$ and $Z$ is the abelian $(\mathbb Z/5\mathbb Z)^8$-covering
$g:Z\to\mathbb P^2$ defined by the field extension $g^*:\mathbb
C(\mathbb P^2)\hookrightarrow \mathbb C(Z)=\mathbb C(\mathbb
P^2)(w_1,\dots,w_8)$ such that $w_j^5=l_j\l_9^{-1}$, for $j=1,\dots,
8$, and $l_j=0$ are equations of $L_j$. The divisors
$\frac{1}{5}h_i^*((\sigma^*(l_jl_9^{-1})))$ belong to
$\text{Tor}_5(X_i)$ and generate in it a subgroup $G_i\simeq
(\mathbb Z/5\mathbb Z)^6$.

From Proposition \ref{2K} it follows that
the space $\mathcal M_{k_{X,d,m},p_{X,d,m}}$ is non-empty, since by
Bombieri Theorem
the surface $X$ being of general type contains smooth curves
numerically equivalent to $dmK_X$ if $dm\geq 5$.

By Proposition 4.1 in \cite{KK1},
$\text{Kl}(X_1)=\text{Aut}(X_1)=Gal(\widetilde X_1/\mathbb P^2)$ and
the action of $\text{Aut}(X_1)$ on $G_1$ is trivial since
$\text{Kl}(X_1)$ leaves fixed the lines $L_j\subset \widetilde L$.
Therefore, by Theorem \ref{main2}, the surfaces $Y_{1,k}$  obtained
as totally ramified numerically $(d,m)$-canonical coverings
$f_k:Y_{1,k}\to X_1$ defined by the divisors $C_k\sim
mK_X+\alpha_k$, $\alpha_k\in G_1$, are not pairwise orientation
preserving diffeomorphic and hence they belong to distinct connected
components of the moduli space $\mathcal M_{k_{X,d,m},p_{X,d,m}}$.

According to Theorem \ref{newTh}, the surfaces $Y_1$ and $Y_2$
obtained as totally ramified numerically $(d,m)$-canonical cyclic
coverings $f_1:Y_1\to X_1$ and $f_2:Y_2\to\bar X_1$ can not be
deformation equivalent, since $\text{Kl}(X_1)=\text{Aut}(X_1)$.
Therefore,
the totally ramified numerically $(d,m)$-canonical cyclic coverings
$f_{1,k}:Y_k\to \overline X_1$ defined by the divisors $C_k\sim
mK_{\overline X_1}+\alpha_k$, $\alpha_k\in G_1$, define another
$5^6$ distinct connected components of the moduli space $\mathcal
M_{k_{X,d,m},p_{X,d,m}}$.

Once more, by Proposition 4.1 in \cite{KK1},
$\text{Aut}(X_2)=Gal(\widetilde X_2/\mathbb P^2)$ and the action of
$\text{Aut}(X_2)$ on $G_2$ is trivial, while
$\text{Out}(\pi_1(X_2))=\text{Kl}(X_2)\ne \text{Aut}(X_2)$. Hence,
the totally ramified numerically $(d,m)$-canonical cyclic coverings
$f:\tilde Y\to X_2$ give at least $5^6$ another distinct connected
components in $\mathcal M_{k_{X,d,m},p_{X,d,m}}$, since these
surfaces $\tilde Y$ are not homeomorphic to the surfaces $Y$
obtained as totally ramified numerically $(d,m)$-canonical cyclic
coverings $f:Y\to X_1$ (the fundamental groups of surfaces $\tilde
Y$ and $Y$ have non isomorphic groups of outer automorphisms). \qed

\begin{thm} \label{main3} Let $f_1 :Y_1\to X$ be a totally ramified numerically
$(2,m)$-canonical cyclic covering of a Miyaoka-Yau  surface $X$ and
let $\widetilde Y_2$ be a surface deformation equivalent to $Y_1$.
Then the canonical model $Y_2$ of $\widetilde Y_2$ can be
represented as a numerically $(2,m)$-canonical totally ramified
cyclic covering of $X$ branched along a curve $B_2\subset X$, where
$B_2$ is a reduced {\rm (}not necessary irreducible{\rm )} curve
with $ADE$-singularities.

If $2m\geq 5$, then the connected component $M$ of the moduli space
to which $Y_1$ belongs is an irreducible variety of dimension
$m(2m-1)K_X^2+p_g(X)$. The Kodaira dimension $\kappa (M)$ of $M$ is
equal to $-\infty$ and if the irregularity $q(X)=0$, then $M$ is an
unirational variety.
\end{thm}
\proof Let $f_1:Y_1\to X$ be defined by a divisor $C\sim
mK_X+\alpha_1$ and branched along a curve $B_1\in |2C|$.

Since $Y_1$ and $\widetilde Y_2$ are deformation equivalent, there
exists an orientation preserving diffeomorphism $\varphi: \widetilde
Y_2\to Y_1$ such that $\varphi^*(K_{Y_1})=K_{\widetilde Y_2}$.
Hence, by Siu rigidity \cite{Siu}, the composition $f_1\circ\varphi$
is homotopic to a holomorphic map $h:\widetilde Y_2\to X$.

We have $\deg h=\deg (f_1\circ\varphi)=2$, and therefore $h^*:
\mathbb C(X)\hookrightarrow \mathbb C(\widetilde Y_2)$ is a Galois
extension of degree $2$. Let $f_2:Y_2\to X$ be the normalization of
$X$ in the field $\mathbb C(\widetilde Y_2)$ and let $\nu
:\widetilde Y_2\to Y_2$ be a morphism such that $h=f_2\circ \nu$ (in
fact, $\nu$ is the minimal resolution of singularities of $Y_2$).

Let us show that $Y_2$ has at most ADE-singularities. Indeed, if
$D\subset \widetilde Y_2$ is an irreducible curve such that $\nu(D)$
is a point, then
$$(K_{\widetilde Y_2},D)_{\widetilde Y_2}=(\varphi^*(K_{Y_1}),D)_{\widetilde Y_2}=
(\varphi^*(f_1^*(D_1)),D)_{\widetilde Y_2}=(h^*(D_1),D)_{\widetilde
Y_2}=0,$$ where $D_1\equiv mK_X+C\equiv (m+1)K_X$. It implies, by
adjunction, that $D$ is a rational curve with $(D^2)_{\widetilde
Y_2}=-2$.

Since $f_2 : Y_2\to X$ is a double covering, its branch curve, which
we denote by $B_2$, has similar to $Y_2$ at most ADE-singularities.
Let us show that $B_2\equiv 2mK_X$. Indeed, we have $2K_{\widetilde
Y_2}=h^*(2K_X+B_2)$. On the other hand,
$$2K_{\widetilde Y_2}= 2\varphi^*(K_{Y_1})\equiv
2(m+1)\varphi^*(f_1^*(K_{X}))\equiv 2(m+1)h^*(K_X)$$ and hence
$B_2\equiv 2mK_X$.

By Remark \ref{rem3}, the covering $f_2: Y_2\to X$ is defined by a
divisor $C_2\sim mK_X+\alpha_2$. A small deformation of $B_2$ into a
smooth curve $B_3\sim B_2$ defines a numerically $(2,m)$-canonical
(totally ramified) covering $f_3:Y_3\to X$ branched along $B_3$ and
associated with the divisor $C_2\sim mK_X+\alpha_2$. Because of the
existence of simultaneous resolution for simple singularities, $Y_2$
and $Y_3$ belong to the same irreducible component of the moduli
space.

By Corollary \ref{cor-main2} and Theorem \ref{main2}, there is an
automorphism $\psi\in \text{Aut}(X)$ such that
$\psi^*(\delta(\alpha_2))=\delta(\alpha_1)$. Therefore, we can
assume that $\delta(\alpha_2)=\delta(\alpha_1)$ (changing $f_2$ by
$\psi\circ f_2$) and furthermore all the surfaces $Y_t$ deformation
equivalent to $Y_1$ are numerically $(2,m)$-canonical totally
ramified cyclic coverings of $X$ defined by $C_t\sim mK_X+\alpha_t$,
where $\alpha_t$ is such that $\delta(\alpha_t)=\delta(\alpha_1)$,
and branched along curves $B_t\sim 2(mK_X+\alpha_t)$ that have at
most ADE-singularities.

Let $\Delta=\delta(\alpha_1)\in\text{Tor}(X)$. Denote by
$\text{Pic}_\Delta(X)=\delta^{-1}(\Delta)\subset \text{Pic}(X)$. As
in the proof of Theorem \ref{def-eq}, consider the scheme
$T_{\Delta, 2m}$ parametrizing the curves $B_t$ in $X$ numerically
equivalent to $2mK_X$ and such that $\delta(B_t)=
2m\delta(K_X)+2\Delta$. The scheme $T_{\Delta,2m}$ is fibered over
$\text{Pic}_{\Delta}(X)$, $\gamma_{1}:
T_{\Delta,2m}\to\text{Pic}_{\Delta}(X)$, with fibers $\gamma_1^{-1}(
\alpha_t)=\mathbb P(H^0(X, \mathcal O_X(2mK_X+ d \alpha_t))$ over $
\alpha_t\in \text{Pic}_{\Delta}(X)$. Obviously, the subscheme $U$
consisting of points $t$ of $T_{\Delta,2m}$ for which $B_t$ are
reduced curves with at most ADE-singularities is a Zariski open
non-empty subset.

Similarly, let $T_{\Delta, 2}$ be the scheme parametrizing the
curves $D_t$ in $X$ numerically equivalent to $2K_X$ and such that
$\delta(D_t)= 2\delta(K_X)+\Delta$. The scheme $T_{\Delta,2}$ also
is fibered over $\text{Pic}_{\Delta}(X)$, $\gamma_{2}:
T_{\Delta,2}\to\text{Pic}_{\Delta}(X)$, with fibers
$\gamma_2^{-1}(\alpha_t)=\mathbb P(H^0(X, \mathcal O_X(2K_X+
\alpha_t))$ over $\widetilde \alpha_t\in \text{Pic}_{\Delta}(X)$.
Denote by $T=U\times_{\text{Pic}_{\Delta}(X)}T_{\Delta,2}$ the
product of the fibrations $\gamma_1$ and $\gamma_2$ and let $p_i$,
$i=1,2$, be the projections of $T$ onto the factors.

Note that $T$ is an irreducible variety and it follows from the
above consideration and the proof of Theorem \ref{def-eq} that the
points of $T$ parametrize all the surfaces deformation equivalent to
$Y_1$. Therefore, the above considerations show that the connected
component $M$ of the moduli space to which $Y_1$ belongs is an
irreducible variety.

To complete the proof of Theorem, consider the surjective morphism
$\mu: T\to M$. It is easy to see that the fibers of $p_1$ are
subvarieties of the fibers of $\mu$. Let us show that each fiber of
$\mu$ is the union of a finite number of fibers of $p_1$. Indeed, a
degree two totally ramified covering $f:Y\to X$ branched along a
curve $B\equiv 2mK_X$ defines a degree two extension $\mathbb
C(X)\hookrightarrow \mathbb C(Y)$. This extension defines an
automorphism $f^*\in \text{Aut}(Y)$ of order two whose set of fixed
points is the ramification locus $R$ of $f$ (in the case when $Y$
has no $(-2)$-curves). Conversely, each $h\in \text{Aut}(Y)$
determines uniquely  the set of its fixed points. But,
$\text{Aut}(X)$ and $\text{Aut}(Y)$ are finite groups, since $X$ and
$Y$ are surfaces of general type. Therefore each fiber of $\mu$ is
the union of a finite number of fibers of $p_2$ and, by Stein
Factorization Theorem, the morphism $\mu$ factorizes through a
finite morphism $\mu_1: U\to M_1$.

We have $\dim M= \dim U =\dim |2mK_X+ d \alpha_t|+q(X)$. By Mumford
Vanishing Theorem and Riemann-Roch Theorem,
$$ \dim |2(mK_X+\alpha_1)|=  \displaystyle
\frac{(2(mK_X+\alpha_t),2 (mK_X+\alpha_t)-K_X)_X}2 +p_g(X)-q(X).$$
Therefore $\dim M= m(2m-1)K^2_X+p_g(X)$. Moreover, since $U$ has
Kodaira dimension  $\kappa(U)=-\infty$, $M$ also has the same
Kodaira dimension. If $q=0$, then $U$ is a rational variety and,
therefore, $M$ is unirational. \qed

\begin{rem} \label{last-rem} {\rm The same
arguments as in the proof of Theorem \ref{main3} show that for any
$d\geq 2$ and for any surface $Y_2$ deformation equivalent to a
surface $Y_1$ that is obtained as a numerically $(d,m)$-canonical
totally ramified cyclic covering of a Miyaoka-Yau surface $X$, there
exists a degree $d$ morphism $f_2:Y_2\to X$; but if $d\geq 3$, then
$f_2$ is not necessary a cyclic covering. Furthermore, these
arguments give rise to a lower bound on the dimension of that
irreducible component $M$ of the moduli space of surfaces to which
$Y_1$ belongs:
$$\displaystyle \dim M\geq
\frac{dm(dm-1)}{2}K_X^2+p_g(X).$$}
\end{rem}

\begin{rem} \label{last-rem} {\rm The proof of Theorem \ref{main3} shows that for any surface $Y$ that is a totally ramified
numerically $(2,m)$-canonical cyclic covering of a Miayoka-Yau
surface $X$, the action of the group $\mathbb Z/2\mathbb Z$ on $Y$
by deck transformation deforms simultaneously with any deformation
of complex structure.}
\end{rem}

\section{New examples of surfaces having no any anti-holomorphic
automorphism} In this section, by the {\it Mostow strong rigidity}
of a compact complex manifold $X$ we mean the following property:
{\it whatever is a homotopy equivalence $p: X\to X$, it is homotopic
to a holomorphic or anti-holomorphic map $X\to X$.} By Mostow
rigidity theorem, the fake projective planes and the surfaces from
Examples 1 and 2 in \cite{KK1} are Mostow strongly rigid. In
addition, they are $K(\pi,1)$ as topological spaces.

Let us underline that here we do not include the unicity statement
into the definition of Mostow strong rigidity; the reason is that
such a unicity statement is not used in the proofs of the results
below. (Note, however, that if the definition had included the
unicity, then in the statement given in Remark \ref{topi} the "if"
could be replaced by the "if and only if".)

\begin{thm}\label{main1}
Let $X$ be a Mostow strongly rigid surface of general type. If $X$
is a $K(\pi,1)$ and $\operatorname{Out}(\pi)$ contains no elements
of even order, then no numerically $(d,m)$-canonical totally
ramified cyclic covering $Y$ of $X$ can be deformed to a surface
isomorphic to $\overline Y$. In particular, all the surfaces $Y'$
obtained by deformation of such a  $Y$ have no any anti-holomorphic
automorphism.
\end{thm}
\proof We argue by contradiction. Assume that $Y$ is deformed to a
surface $Y'$ that is isomorphic to $\overline Y$. Identify $Y'$ as a
smooth manifold with $Y$ (in accordance with the deformation between
them) and denote by $c: Y\to Y'$ an anti-holomorphic diffeomorphism
between $Y$ and $Y'$. Preserve the same notation $c$ for the
diffeomorphism induced by $c$ on $Y$ following our identification of
the smooth manifolds underlying $Y$ and $Y'$.

Consider  the element $c_*\in\operatorname{Out} (\pi_Y)$ induced by
$c$ (and the deformation equivalence between $Y$ and $Y'$). By
Proposition  \ref{pi}, $\pi_1(Y)= \pi_1(X)$. Since $X$ is Mostow
strongly rigid, is of general type, and is $K(\pi,1)$, the element
$c_*$,  as any element of $\operatorname{Out}
(\pi_X)=\operatorname{Out} (\pi_Y)$, is of finite order. Denote by
$n$ the order of $c_*$, $n>0$. Since $X$ is a $K(\pi,1)$,
Proposition \ref{pi} implies that the maps $f: Y'=Y\to X$ and
$f\circ c^n$ are homotopic, and since in addition
$c^*(\delta(K_{Y}))=-\delta(K_{Y'})=-\delta(K_{Y})$ (the Chern class
$\delta(K)$ being integral does not change under deformations) while
$f^*$ transforms $(dm-m+1)\delta(K_X)$ in $\delta(K_Y)$ modulo
torsion, it implies that $n$ is even. Contradiction.\qed

\begin{rem}\label{symplectic}
{\rm Literally the same proof shows that the surfaces $Y$ (and $Y'$)
satisfying Theorem \ref{main1} assumptions have no diffeomorphisms
$f:Y\to Y$ with $f^*[K]=-[K], [K]\in H^2(Y;\mathbb Q)$. Hence, due
to the invariance of the canonical class under deformations in the
class of almost-complex structures, the imaginary part $\omega$ of
the K\"ahler structure of $Y$ considered as a symplectic structure
on the underlying smooth manifold is not symplectic deformation
equivalent to its reverse, $-\omega$. Thus, Theorem \ref{main1} and,
in particular, its Corollary \ref{main-ex} below provide new
examples of opposite symplectic structures not equivalent to each
other (cf., \cite{KK2}).}
\end{rem}

\begin{rem}\label{topi}
{\rm If $X$ is a Mostow strongly rigid surface of general type and
$X$ is a $K(\pi, 1)$, then the fundamental group $\pi_1(X)$ has no
automorphisms of even order $>0$ as soon as $X$ has no neither any
anti-holomorphic automorphism or any holomorphic automorphism of
non-zero even period. }
\end{rem}

\begin{cor} \label{main-ex} If $X$ is a fake projective plane or the rigid surface constructed
in Example 1 in \cite{KK1}, then each of the surfaces $Y'$ obtained
by deformation of a numerically $(d,m)$-canonical totally ramified
cyclic covering $Y$ of $X$ can not be deformed to its complex
conjugate, and, in particular, has no any anti-holomorphic
automorphism.
\end{cor}
\proof As is proved in \cite{KK1}, each of these surfaces $X$ has
neither anti-holomorphic automorphisms nor holomorphic automorphisms
of even order $>0$, so that the result follows from Theorem
\ref{main1} and Remark \ref{topi}. \qed

\ifx\undefined\bysame
\newcommand{\bysame}{\leavevmode\hbox to3em{\hrulefill}\,}
\fi

\end{document}